\documentclass{amsart}
\usepackage{amssymb, amsmath}
\usepackage[dvips]{graphicx}
\usepackage{amsfonts}
\usepackage{latexsym}
\usepackage{color}
\newtheorem{theorem}{Theorem}	
\newtheorem{lemma}{Lemma}[section]		
		
\newtheorem{proposition}{Proposition}		
			
\newtheorem{definition}{Definition}	
\newtheorem{example}{Example}
\newtheorem{problem}{Problem}

\begin{document}

\title{Self-dual Wulff shapes and spherical convex bodies 
of constant width ${\pi}/{2}$
}
\thanks{\color{black} This work was partially supported by JSPS KAKENHI Grant Number 26610035.}
\author{
Huhe Han}
\address{Graduate School of Environment and Information Sciences,Yokohama National University, {Yokohama 240-8501,} Japan}
\email{han-huhe-bx@ynu.jp}
\author{Takashi Nishimura
}
\address{
Research Institute of Environment and Information Sciences,  
Yokohama National University, 
Yokohama 240-8501, Japan}
\email{nishimura-takashi-yx@ynu.jp}



\begin{abstract}
For any Wulff shape, 
its dual Wulff shape 
is naturally defined.   
A self-dual Wulff shape is a Wulff shape equaling its dual Wulff shape exactly.   
In this paper, it is shown that a Wulff shape 
is self-dual if and only if the spherical convex body induced by it 
is of constant width 
${\pi}/{2}$.
\end{abstract}
\maketitle

\section{Introduction}
For a positive integer $n$, let $S^n$ be the unit sphere in $\mathbb{R}^{n+1}$.   
Let $\mathbb{R}_+$ be the set consisting of positive real numbers.
For any continuous function $\gamma: S^n\to \mathbb{R}_+$ and 
any $\theta\in S^n$, 
let $\Gamma_{\gamma, \theta}$ be the set consisting of 
$x\in \mathbb{R}^{n+1}$ such that 
$x\cdot \theta\le \gamma(\theta)$, 
where the dot in the center stands for the scalar product of two vectors 
$x, \theta\in \mathbb{R}^{n+1}$.    
Then, the \textit{Wulff shape} associated with 
the support function $\gamma$ is 
the following set $\mathcal{W}_\gamma$:  
\[
\mathcal{W}_\gamma=\bigcap_{\theta\in S^n}\Gamma_{\gamma, \theta}.
\] 
A Wulff shape $\mathcal{W}_\gamma$ was firstly introduced by G. Wulff in \cite{wulff} 
as a geometric model of a crystal at equilibrium.  
By definition, any Wulff shape is a convex body in $\mathbb{R}^{n+1}$ containing the origin as an interior point.   
Conversely, it has been known that for any convex body $W$ in $\mathbb{R}^{n+1}$ such that int$(W)$ contains the origin where 
int$(W)$ stands for the set of interior points of $W$, there exists a continuous function $\gamma: S^n\to \mathbb{R}_+$ 
such that $W=\mathcal{W}_\gamma$ (\cite{taylor}).     
By using the polar plot expression of elements of 
$\mathbb{R}^{n+1}-\{0\}$, $S^n\times \mathbb{R}_+$ 
may be naturally identified with $\mathbb{R}^{n+1}-\{0\}$.    
Under this identification, 
for any Wulff shape $\mathcal{W}_\gamma$ and any $\theta\in S^n$, 
the intersection $\partial \mathcal{W}_\gamma \cap L_\theta$ 
is exactly one point (denoted by $(\theta, w(\theta))$), 
where $\partial \mathcal{W}_\gamma$ is the boundary of 
$\mathcal{W}_\gamma$ and $L_\theta$ is the half line 
$L_\theta=\{(\theta, r)\; |\; r\in \mathbb{R}_+\}$.
For a Wulff shape $\mathcal{W}_\gamma$, 
let $\overline{\gamma}: S^n\to \mathbb{R}_+$ be the continuous function 
defined by $\overline{\gamma}(\theta)=\frac{1}{w(-\theta)}$.     
Then, the Wulff shape $\mathcal{W}_{\overline{\gamma}}$ is called 
the {\it dual Wulff shape} of $\mathcal{W}_\gamma$ and is denoted by 
$\mathcal{DW}_\gamma$.   
For any Wulff shape $\mathcal{W}_\gamma$, 
there is a characterization of the dual Wulff shape of $\mathcal{W}_\gamma$.    
The graph of a continous function 
$\gamma: S^n\to \mathbb{R}_+$ is denoted by 
graph$(\gamma)$.    Let 
inv$: \mathbb{R}^{n+1}-\{0\}\to \mathbb{R}^{n+1}$ be the inversion of 
$\mathbb{R}^{n+1}-\{0\}$ defined by inv$(\theta, r)=(-\theta, \frac{1}{r})$.    
Then, for any continuous function 
$\gamma: S^n\to \mathbb{R}_+$, 
$\mathcal{DW}_\gamma$ is exactly 
the convex hull of inv(graph$(\gamma)$). 
By this characterization,  it is clear that 
$\mathcal{DDW}_\gamma$ is $\mathcal{W}_\gamma$ 
for any $\mathcal{W}_\gamma$ when inv(graph$(\gamma)$) is the boundary of the convex hull of inv(graph$(\gamma)$).           
A Wulff shape $\mathcal{W}_\gamma$ is said to be {\it self-dual} if the equality 
$\mathcal{W}_\gamma=\mathcal{DW}_\gamma$ holds.   
\par 
\medskip 
In this paper, a simple and useful characterization for a self-dual Wulff shape in $\mathbb{R}^{n+1}$ is given.    
In order to state our characterization, 
several notions in $S^{n+1}$ are defined.  
For any point $P$ of $S^{n+1}$, 
let $H(P)$ be the hemisphere centered at $P$, namely $H(P)$ is the subset 
of $S^{n+1}$ consisting of $Q\in S^{n+1}$ satisfying $P\cdot Q\ge 0$, where 
the dot in the center stands for the scalar product of 
two vectors $P, Q\in \mathbb{R}^{n+2}$.  
A subset $\widetilde{W}$ of $S^{n+1}$ is said to be {\it hemispherical} if there 
exists a point $P\in S^{n+1}$ such that 
$\widetilde{W}\cap H(P)=\emptyset$.    
A hemispherical subset $\widetilde{W}\subset S^{n+1}$ is said to be 
{\it spherical convex} if for any $P, Q\in \widetilde{W}$ the following arc $PQ$ 
is contained in $\widetilde{W}$: 
\[
PQ=
\left\{
\left. 
\frac{(1-t)P+tQ}{||(1-t)P+tQ||}\; \right|\; t\in [0,1]
\right\}.   
\] 
A hemispherical subset $\widetilde{W}$ is called a {\it spherical convex body} 
if it is closed, spherical convex and has an interior point.    
A hemisphere $H(P)$ is said to {\it support a spherical convex body $\widetilde{W}$} if both $\widetilde{W}\subset H(P)$ 
and $\partial \widetilde{W}\cap \partial H(P)\ne \emptyset$ hold.
For a spherical convex body $\widetilde{W}$ and a hemisphere $H(P)$ supporting $\widetilde{W}$, 
following  \cite{lassak, lassak2}, the {width of 
$\widetilde{W}$ determined by $H(P)$} 
is defined as follows.    
For any two $P,Q\in S^{n+1}$ $(P\ne \pm Q)$, the intersection 
$H(P)\cap H(Q)$ is called a {\it lune} of $S^{n+1}$.    
The {\it thickness of the lune $H(P)\cap H(Q)$}, 
denoted by $\triangle(H(P)\cap H(Q))$, 
is the real number 
$\pi - |PQ|$, where $|PQ|$ stands for the length of the arc $PQ$.   
For a spherical convex body 
$\widetilde{W}$ and a hemisphere $H(P)$ supporting $\widetilde{W}$, 
the {\it width of $\widetilde{W}$ determined by $H(P)$}, denoted by 
$\mbox{width}_{H(P)}\widetilde{W}$, is the minimum of the following set:   
\[
\left\{
\triangle(H(P)\cap H(Q))\; \left|\; \widetilde{W}\subset H(P)\cap H(Q), 
H(Q) \mbox{ supports }\widetilde{W}
\right.\right\}.   
\]  
For any $\rho\in \mathbb{R}_+$ less than $\pi$,   
a spherical convex body $\widetilde{W}\subset S^{n+1}$ is said to be 
{\it of constant width $\rho$} if 
$\mbox{width}_{H(P)}\widetilde{W}
=\rho$ for any $H(P)$ supporting $\widetilde{W}$.   
\par 
\medskip 
Let $Id: \mathbb{R}^{n+1}\to \mathbb{R}^{n+1}\times 
\{1\}\subset \mathbb{R}^{n+2}$, $N\in S^{n+1}$ and  
$\alpha_N: S^{n+1}-H(-N)\to 
\mathbb{R}^{n+1}\times \{1\}\subset \mathbb{R}^{n+2}$ 
  be 
the mapping defined by $Id(x)=(x,1)$, the point 
$(0, \ldots, 0,1)\in S^{n+1}$ and  
the central projection defined as follows respectively.   
\begin{eqnarray*}
{ } & { } \alpha_N\left(P_1, \ldots, P_{n+1}, P_{n+2}\right)
=
\left(\frac{P_1}{P_{n+2}}, \ldots, \frac{P_{n+1}}{P_{n+2}}, 
{\color{black}1}\right)\\ 
{ } & { } \qquad \qquad \qquad 
(\forall \left(P_1, \ldots, P_{n+1}, P_{n+2}\right)\in S^{n+1}-H(-N)). 
\end{eqnarray*} 
Then, for any Wulff shape $\mathcal{W}_\gamma$, it is clear that 
$\alpha_N^{-1}\circ Id(\mathcal{W}_\gamma)$ is a spherical convex body.  
The set $\alpha_N^{-1}\circ Id(\mathcal{W}_\gamma)$ is called the 
{\it spherical convex body induced by} $\mathcal{W}_\gamma$.    
\begin{theorem}\label{theorem 1}
Let $\gamma: S^n\to \mathbb{R}_+$ be a continuous function.   
Then, the Wulff shape $W_\gamma$ is self-dual if and only if 
the spherical convex body induced by $W_\gamma$ 
is of constant width ${\pi}/{2}$.
 \end{theorem}
\noindent 
The unit disc $D^{n+1}=\{x\in \mathbb{R}^{n+1}\;  |\; ||x||\le 1\}$ of 
$\mathbb{R}^{n+1}$ is clearly self-dual.    
Let $R$ be a rotation of $\mathbb{R}^{n+2}$ about an $n$ dimensional 
linear subspace with a small angle.   
Then, since the property of constant width is an invariant property 
by $R$, 
by Theorem \ref{theorem 1}, 
$Id^{-1}\circ \alpha_N
\left(R\left(\alpha_N^{-1}\circ Id(D^{n+1})\right)\right)$ 
is self-dual as well 
(see Figure 1).  
\begin{figure}
\begin{center}
\includegraphics[width=10cm]{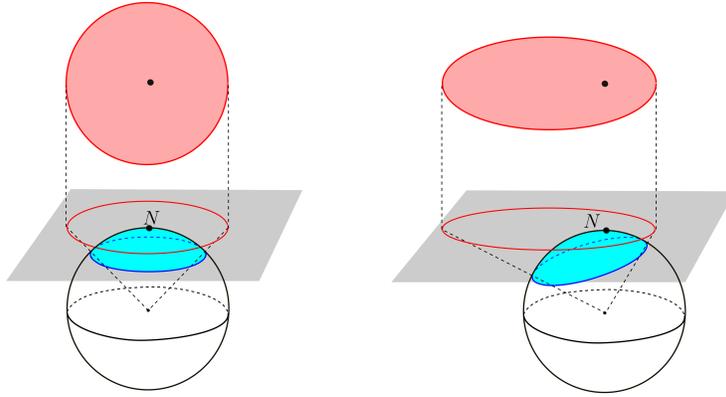}
\caption{
{\color{black}Self-dual Wulff shapes include central projections of spherical caps of 
width $\pi/2$.}   
}
\label{figure1}
\end{center}
\end{figure}
Moreover, 
let $\widetilde{\triangle}$ be a spherical 
triangle of constant width 
$\frac{\pi}{2}$ in $S^2$ containing $N$ as an interior point.     
Then, by Theorem \ref{theorem 1}, 
not only $Id^{-1}\circ \alpha_N\left(\widetilde{\triangle}\right)$ itself, 
but also 
any $Id^{-1}\circ \alpha_N
\left(R\left(\widetilde{\triangle}\right)\right)$ 
is self-dual (see Figure 2).   
{\color{black}For more consideration on simple, explicit examples,  
see Section \ref{section 4}.}
\begin{figure}
\begin{center}
\includegraphics[width=10cm]{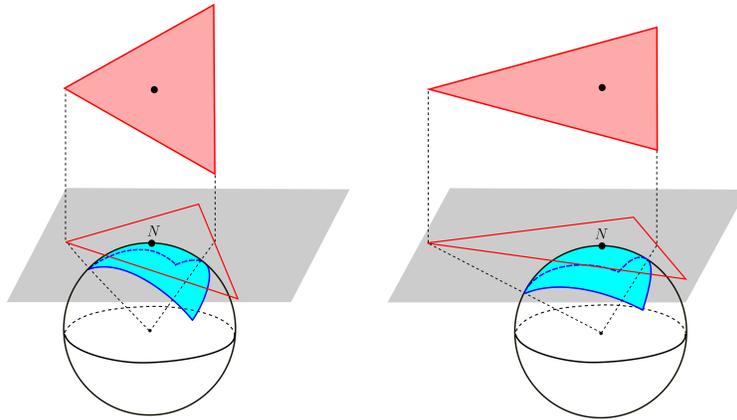}
\caption{
{\color{black}
Self-dual Wulff shapes include triangles which are central projections of constant-width 
spherical triangles of width $\pi/2$.}
}
\label{figure2}
\end{center}
\end{figure}
\par 
On the other hand, any Reuleaux triangle in $\mathbb{R}^2$ containing the origin as an interior point (see Figure 3) is not a self-dual Wulff shape, 
although it is a Wulff shape of constant width {\color{black}in $\mathbb{R}^2$}.   
This is because any Reuleaux triangle is strictly convex, 
and thus the boundary of it must be smooth by \cite{hannishimura2} 
if it is self-dual.   
However, there are three non-smooth points for any Reuleaux triangle in $\mathbb{R}^2$.  
{\color{black} By Theorem 1, its spherical convex body is not of constant width $\pi/2$}.   
\begin{figure}
\begin{center}
\includegraphics[width=4cm]{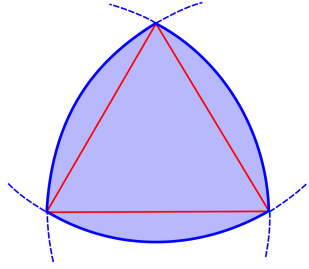}
\caption{
Reuleaux triangle. 
}
\label{figure3}
\end{center}
\end{figure}
\par 
\bigskip  
  In Section \ref{section 2}, preliminaries for the proof of Theorem 
  \ref{theorem 1} are given.  
The proof of Theorem \ref{theorem 1} is given in Section \ref{section 3}.     
{\color{black}Finally, more consideration on simple, explicit examples 
is given.} 
\section{Preliminaries}\label{section 2}
The follwing two theorems given in \cite{lassak} are keys 
for the proof of Theorem \ref{theorem 1}.
\begin{theorem}[\cite{lassak}]\label{center}
Let $\widetilde{W}\subset S^{n+1}$ be a spherical convex body and 
let H(P) be a hemisphere which supports 
$\widetilde{W}$.
\begin{enumerate}
\item If $P\notin \widetilde{W}$, then there exists a unique hemisphere 
$H(Q)$ supporting $\widetilde{W}$ 
such that the lune $H(P)\cap H(Q)$ contains $\widetilde{W}$ 
and has thickness {\rm width}$_{P}(\widetilde{W})$. 
This hemisphere supports $\widetilde{W}$ at the point 
$R$ at which the largest ball $B(P, r)$ touches $\widetilde{W}$ from outsite.
We have $\Delta (H(P)\cap H(Q))=\frac{\pi}{2}-r$.
\item If $P\in \partial \widetilde{W}$, then there exists at least one hemisphere
$H(Q)$ supporting $\widetilde{W}$ such that $H(P)\cap H(Q)$ is a lune containing 
$\widetilde{W}$ which has thickness {\rm width}$_{P}(\widetilde{W})$. 
This hemisphere supports $\widetilde{W}$ at $R = P$. 
We have $\Delta (H(P)\cap H(Q))=\frac{\pi}{2}$.
\item If $P\in {\rm int}(\widetilde{W})$, then there exists at least one hemisphere $H(Q)$ supporting
$\widetilde{W}$ such that $H(P)\cap H(Q)$ is a lune containing $\widetilde{W}$ which has thickness 
{\rm width}$_{P}(\widetilde{W})$. Every such $H(Q)$ supports $\widetilde{W}$ at exactly one point
$R\in \partial \widetilde{W}\cap B(P, r)$, where B(P, r) dentes the largest ball with center $P$ 
contaned in $\widetilde{W}$, and for every such $R$ this hemisphere $H(Q)$,
denoted $H_{R}(Q)$, is unique. For every $R$ we have 
$\Delta (H(P)\cap H(Q))=\frac{\pi}{2}+r$.
\end{enumerate}
\end{theorem}
\begin{definition}[\cite{lassak}]\label{diameter}
{\rm 
Let $\widetilde{W}\subset S^{n+1}$ be a spherical convex body.   
Then, the following real number is called the {\it diameter} of $\widetilde{W}$ and is denoted by 
diam$\left(\widetilde{W}\right)$.  
\[
\max\left\{|PQ|\; \left|\; P,Q\in \widetilde{W}\right.\right\}.   
\] 
}
\end{definition}
\begin{theorem}[\cite{lassak}]\label{max}
Let $\widetilde{W}\subset S^{n+1}$ be a spherical convex body.   
Suppose that {\rm diam} $\left(\widetilde {W}\right)\leq \frac{\pi}{2}$.   
Then, 
the following holds:   \\
\indent ${\rm diam}\left(\widetilde{W}\right)=\\$
\indent \ \ \ ${\rm max} \left\{\left.{\rm width}_{H(P)}\left(\widetilde{W}\right) \; \right|\; H(P)\ {\rm is\ a\ supporting\ hemisphere\ of\ } \widetilde{W} \right\}$.
\end{theorem}

\par
\begin{definition}[\cite{nishimurasakemi2}]\label{definition 2.3}
{\rm 
For any hemispherical subset $\widetilde{W}$ of $S^{n+1}$, 
the following set (denoted by $\mbox{\rm s-conv}\left(\widetilde{W}\right)$) is called the {\it spherical convex hull of} 
$\widetilde{W}$:     
\[
\mbox{\rm s-conv}\left(\widetilde{W}\right)= 
\left\{\left.
\frac{\sum_{i=1}^k t_iP_i}{||\sum_{i=1}^kt_iP_i||}\;\right|\; 
P_i\in \widetilde{W},\; \sum_{i=1}^kt_i=1,\; t_i\ge 0, k\in \mathbb{N}
\right\}.
\]  
}
\end{definition}
It is clear that 
$\mbox{\rm s-conv}\left(\widetilde{W}\right)=\widetilde{W}$ if 
$\widetilde{W}$ is spherical convex.      More generally, we have the following:   
\begin{lemma}[\cite{nishimurasakemi2}]\label{lemma 2.4}
Let $\widetilde{W}$ be a hemispherical subset of $S^{n+1}$.   
Then, the spherical convex hull of $\widetilde{W}$ is the smallest spherical convex set containing 
$\widetilde{W}$.
\end{lemma} 
\begin{definition}[\cite{nishimurasakemi2}]\label{definition 1.1}
{\rm 
For any subset $\widetilde{W}$ of $S^{n+1}$, the set 
$$\bigcap_{P\in \widetilde{W}}H(P)$$ 
is called the {\it spherical polar set} of $\widetilde{W}$ and is denoted by $\widetilde{W}^\circ$.   
}
\end{definition}  
\noindent 
For the spherical polar sets, the following lemma is fundamental.    
\begin{lemma}[\cite{nishimurasakemi2}]\label{lemma 1}
For any non-empty closed hemispherical subset $\widetilde{W} \subset S^{n+1}$, 
{\color{black}
the equality $\mbox{ \rm s-conv}\left(\widetilde{W}\right)= \left(
\mbox{ \rm s-conv}\left( \widetilde{W} \right) \right)^{\circ\circ }$ holds.
} 
\end{lemma}
%
\section{Proof of Theorem \ref{theorem 1}} \label{section 3}
By  the definition of the dual Wulff shape $\mathcal{DW}_\gamma$ for a given Wulff shape $\mathcal{W}_\gamma$, 
it is sufficient to show the following:   
\begin{proposition}\label{constant}
Set $\widetilde{W}=\alpha_N^{-1}\circ Id(\mathcal{W}_\gamma)$.   
Then, $\widetilde{W}=\widetilde{W}^{\circ}$
 if and only if $\widetilde{W}$ is of constant width ${\pi}/{2}$. 
\end{proposition}
\subsection{Proof of the \lq\lq if\rq\rq\;  part of Proposition \ref{constant}}
\quad 
\rm{ 
In this subsection, we show that $\widetilde{W}=\widetilde{W}^{\circ}$
under the assumption that 
$\widetilde{W}$ is of constant width $\frac{\pi}{2}$.
We first show the  inclusion $\widetilde{W}\subset \widetilde{W}^{\circ}$. 
Let $P_1, Q_1$ be two points of $\partial \widetilde{W}$ such that 
$| P_1Q_1 |={\rm diam}(\widetilde {W})$.
Set $P_1=(r\theta, x_{n+2})$ $(0<r, x_{n+2}<1, \theta\in S^n)$.   
Since $\widetilde{W}$ is a spherical convex body, 
for the $\theta\in S^n$, 
there exists the unique real number $t$ $(0<t<1)$ 
such that $H\left(\frac{t\theta+(1-t)N}{||t\theta+(1-t)N||}\right)$ 
supports $\widetilde{W}$.    
For the $t$, set $P=\frac{t\theta+(1-t)N}{||t\theta+(1-t)N||}$.    
%
Then, since we have assumed that 
$\widetilde{W}$ is of constant width $\frac{\pi}{2}$, 
by Theorem \ref{center}, we have that 
${P}\in \partial \widetilde{W}$. 
This implies $P_1={P}$ and hemisphere $H(P_1)$ supports $\widetilde{W}$.
Since $Q_1\in \widetilde{W}\subset H(P_1)$, we have the following,
\[
{\rm diam}(\widetilde{W})=\mid P_1Q_1\mid \leq \frac{\pi}{2}.
 \]  
Let $R$ be an arbitrary point of $\widetilde{W}$. 
Since diam($\widetilde{W}$) $\leq \frac{\pi}{2}$, the following holds,
\[
R\in \bigcap_{\widetilde{R}
\in \widetilde{W}}H(\widetilde{R})=\widetilde{W}^{\circ}.
\]
Therefore, we have $\widetilde{W}\subset \widetilde{W}^{\circ}$.
\par
Next we show the converse inclusion 
$\widetilde{W}^{\circ}\subset \widetilde{W}$.
Suppose that there exists a point $P\in \widetilde{W}^{\circ}$ 
such that $P\notin \widetilde{W}$. 
By Lemma \ref{lemma 1}, it follows that 
$P\notin \widetilde{W}=\cap_{Q\in \widetilde {W}^\circ}H(Q)$.
This implies that 
there exist two points $P$ and $Q$ of $\widetilde{W}^\circ$ 
such that $|PQ|>\frac{\pi}{2}$. 
For these two points $P, Q\in \widetilde{W}^\circ$,  
set $\widetilde{P}=PQ\cap \partial H(P), 
\widetilde{Q}=PQ\cap \partial H(Q)$ (see Figure \ref{figure4}). 
Then we have the following,
\[ 
\pi=
\mid P\widetilde{P}\mid 
+ \mid \widetilde{Q}Q\mid=
\mid P\widetilde{Q}\mid +
\mid \widetilde{Q}\widetilde{P}\mid +
\mid \widetilde{Q}\widetilde{P}\mid+
\mid \widetilde{P} Q\mid =
\mid PQ\mid + \mid \widetilde{P}\widetilde{Q}\mid.
\]
\begin{figure}
\begin{center}
\includegraphics[width=4cm]{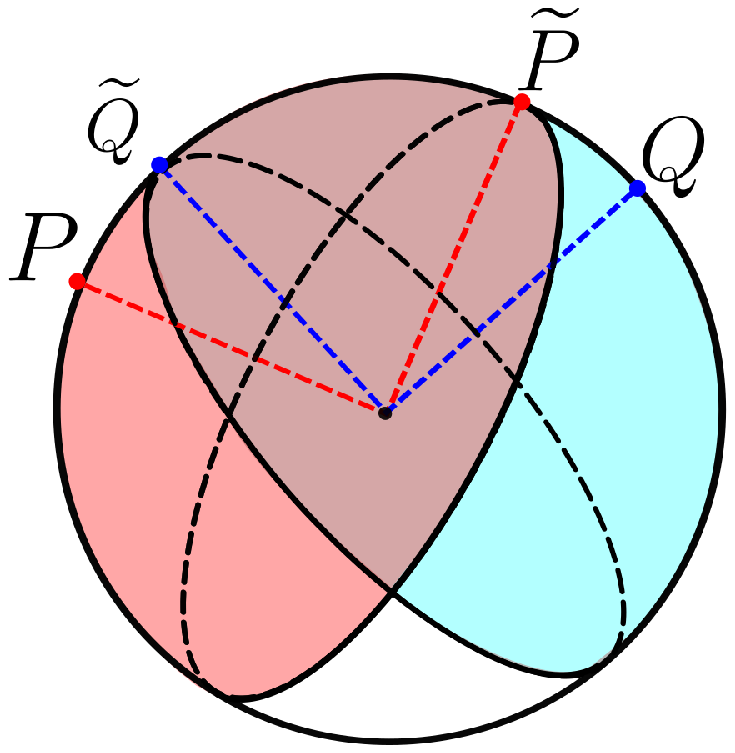}
\caption{
 $\mid PQ\mid >\frac{\pi}{2}$.
}
\label{figure4}
\end{center}
\end{figure}
By the asumption, it follows that
$\mid \widetilde{P}\widetilde{Q}\mid<\frac{\pi}{2}$. 
Let $H(\widetilde{R})$ be 
a suppoting hemisphere of $\widetilde{W}$ whose boundary is  
perpendicular to the arc $PQ$ at the intersecting point.     
Then, the following holds:  
\[
{\rm width}_{H(\widetilde{R})}(\widetilde{W})\leq \mid \widetilde{P}\widetilde{Q}\mid<\frac{\pi}{2}.
\]
This contradicts the assumption that $\widetilde{W}$ is of constant width $\frac{\pi}{2}$.
Therefore, it follows that $\widetilde{W}^{\circ}\subset \widetilde{W}$.
\hfill $\square$\\
\subsection{Proof of the ''only if'' part of Proposition \ref{constant}}
\quad 
\rm{ In this subsection, we show that $\widetilde{W}$ is 
of constant width $\frac{\pi}{2}$ 
under the assumption that 
$\widetilde{W}=\widetilde{W}^{\circ}$. 
Suppose that there exists a hemisphere $H(P)$ supporting $\widetilde{W}$
such that width$_{H(P)}(\widetilde{W})>\frac{\pi}{2}$. 
By Theorem \ref{max}, it follows that diam$(\widetilde{W})\geq {\rm width}_{H(P)}(\widetilde{W})>\frac{\pi}{2}$. 
This implies that there exist two points 
$P, Q \in \widetilde{W}$ sucht that $P\notin H(Q)$. 
Then, we have the following: 
\[
P\notin \bigcap_{Q\in \widetilde{W}}H(Q)=\widetilde{W}^{\circ}.
 \] 
This contradicts the assumption $\widetilde{W}=\widetilde{W}^{\circ}$.
\par
Suppose that there exists a hemisphere $H(P)$ 
supporting $\widetilde{W}$ such that the following holds:  
\[
\mbox{width}_{H(P)}(\widetilde{W})<\frac{\pi}{2}.
\]
Then, there exists a hemisphere $H(Q)$ 
supporting $\widetilde{W}$ such that the following holds:  
\[
\Delta(H(P)\cap H(Q))
={\rm width}_{H(P)}(\widetilde{W}) < \frac{\pi}{2}.
\]
Since 
the thickness $\Delta (H(P)\cap H(Q))=\pi - \mid PQ \mid$, 
we have the following:   
\[ 
\mid PQ \mid> \pi - \frac{\pi}{2}=\frac{\pi}{2}.
\]
On the other hand, 
since $\widetilde{W}$ is a subset of $H(P)$ $(\mbox{\rm resp. } H(Q))$, 
it follows that 
$P\in \widetilde{W}^{\circ}=
\widetilde{W}$ $(\mbox{\rm resp. } Q\in \widetilde{W}^{\circ}=\widetilde{W})$. 
This implies diam($\widetilde{W}$) $\geq \mid PQ \mid>\frac{\pi}{2}$. 
Thus, we have a contradiction.
    \hfill $\square$
\par     
\bigskip 
\noindent 
{\bf Remark.}\quad 
By the proof of Proposition \ref{constant}, it can be seen that 
Proposition \ref{constant} holds for any spherical convex body 
$\widetilde{W}$.     
Namely, we have the following:   
\begin{proposition}\label{self-polar}
Let $\widetilde{W}\subset S^{n+1}$ be a spherical convex body.   
Then, $\widetilde{W}=\widetilde{W}^{\circ}$ 
 if and only if \;$\widetilde{W}$ is of constant width ${\pi}/{2}$. 

\end{proposition}   
{\color{black} 
\section{More on simple, explicit examples}\label{section 4}
\subsection{Centrally symmetric self-dual Wulff shapes}
\label{subsection 4.1}
In this subsection, we determine centrally symmetric Wulff shapes.   
Here, a convex body $W\subset \mathbb{R}^{n+1}$ 
is said to be {\it centrally symmetric} if 
$x\in W$ implies $-x\in W$.    
\begin{proposition}\label{proposition centrally symmetric}
Let $W\subset \mathbb{R}^{n+1}$ be a self-dual Wulff shape. 
Then, $W$ is centrally symmetric if and only if 
$W$ is the unit disc $D^{n+1}$.
\end{proposition}
\indent
{\it Proof}: \quad 
The \lq\lq if\rq\rq part is clear.   
We show the \lq\lq only if\rq\rq part.    
Suppose that there exists a centrally symmetric self-dual Wulff shape $W$ 
which is not the unit disc $D^{n+1}$.   
Then one of the following holds.
\begin{enumerate}
\item[(1)] There exist a point $p\in W$ such that $||p||>1$.
\item[(2)]  The inequality $||p||\leq1$ holds for any point $p$ of $W$  
and there exists a point $q\in \partial W$ such that $||q||<1$.   
\end{enumerate}
\par
\noindent
Here, $||x||$ is the distance from the origin to the point $x\in \mathbb{R}^{n+1}$.
\par
\indent
Suposse that (1) holds.   
Then, since $W$ is centrally symmetric, it follows that $-p\in W$.
Set ${\widetilde{p}=\frac{p}{||p||}\in S^{n}}$.   
For any point $x\in \mathbb{R}^{n+1}$, set 
$X_{+}=\alpha_{N}^{-1}\circ Id(x)$ and $X_{-}=\alpha_{N}^{-1}\circ Id(-x)$
Notice that $P_{-}\in \widetilde{W}=\alpha_N^{-1}\circ Id(W)$.    
Since the distance $|\widetilde{P}_{+}{\widetilde{P}}_{-}|$ is equal to $\frac{\pi}{2}$,
we have the following:
\[
\frac{\pi}{2}=|{\widetilde{P}}_{+}{\widetilde{P}}_{-}|
<|{P}_{+}{P}_{-}|.
\]
This implies $P_{+}\notin H({P}_{-})$.
Thus, it follows that 
\[
{P}_{+}\notin \bigcap_{{Q}\in \widetilde{W}}H({Q})
=\widetilde{W}^{\circ}.
\]
On the other hands, since $W$ is a self-dual Wulff shape and $p\in W$, 
we have that 
${P}_{+}\in \widetilde{W}=\widetilde{W}^{\circ}$. Therefore, we have a contradiction.
\par
\indent
Next, suppose that (2) holds.   
Since there exists a point $q\in \partial W$ such that $||q||<1$, 
it follows that the point $\frac{q}{||q||}\in S^{n}$ does not belong to $W$. 
Set $\widetilde{q}=\frac{q}{||q||}$.     
Then, since $W$ is a self-dual Wulff shape, 
it follows that $\widetilde{Q}_{+}\notin \widetilde{W}=\widetilde{W}^{\circ}$.
On the other hands, by the assumption (2), the following holds.
\[
\widetilde{W}\subset \alpha_N^{-1}\circ Id(D^{n+1})\subset H(\widetilde{Q}_{+}). 
\]
Thus, $\widetilde{Q}_{+}$ is a point of 
$\widetilde{W}^{\circ}$ and we have a contradiction. \hfill $\square$
\subsection{Self-dual Wulff shapes of polytope type}   
A Wulff shape is said to be {\it of polytope type} 
if there exist finitely many points $P_1, \ldots, P_{k}\in S^{n+1}$ such that 
$\widetilde{W}=\bigcap_{i=1}^{k}H(P_i)$, 
where $\widetilde{W}$ is a spherical convex body induced by $W$ and $k\ge n+2\in \mathbb{N}$.    
For crystallines, we have the following proposision: 
\begin{proposition}\label{proposition 4}
Let $W\subset \mathbb{R}^{n+1}$ be a Wulff shape of polytope type   
and let $\widetilde{W}$ be a spherical convex body induced by $W$. 
Set $\widetilde{W}=\bigcap_{i=1}^{k}H(P_i)
\subset S^{n+1}$.  
Then, $W$ is a self-dual Wulff shape if and only if
 $P_{i}$ is a vertex of $\widetilde{W}$ for any $i$ $(1\leq i \leq k)$.
\end{proposition}
\indent
{\it Proof}: \quad 
For the proof, the following lemma is needed.   
\begin{lemma}[Maehara's Lemma \cite{nishimurasakemi2}]
\label{maehara's lemma}
For any hemispherical finite subset $X=\{P_{1}, \ldots , P_{k}\}$, the following holds:
\[
\left\{\left.\frac{\sum_{i=1}^{k}t_{i}P_{i}}{||\sum_{i=1}^{k}t_{i}P_{i}||}\; 
\right|\; P_{i}\in X, \sum_{i=1}^{k}t_{i}=1, t_{i}\geq 0\right\}^{\circ}= 
\bigcap_{i=1}^{k}H(P_{i}).
\]
\end{lemma}
\indent
{\it Proof of  the \lq\lq only if\rq\rq\;  part}: \quad 
Let $W$ be a self-dual Wulff shape of polytope type.    
Then, by {\rm Maehara's} Lemma, we have the following equality: 
\[ 
\widetilde{W}=\bigcap_{i=1}^{k}H(P_{i})=
\left\{\left.\frac{\sum_{i=1}^{k}t_{i}P_{i}}{||\sum_{i=1}^{k}t_{i}P_{i}||}\; 
\right|\; P_{i}\in X, \sum_{i=1}^{k}t_{i}=1, t_{i}\geq 0\right\}^{\circ}.
\]
Then, by Lemma \ref{lemma 1}, the following holds:   
\[ 
\widetilde{W}^\circ=
\left\{\left.\frac{\sum_{i=1}^{k}t_{i}P_{i}}{||\sum_{i=1}^{k}t_{i}P_{i}||}\; 
\right|\; P_{i}\in X, \sum_{i=1}^{k}t_{i}=1, t_{i}\geq 0\right\}.
\]
Since $W$ is a self-dual Wulff shape, it follows that 
$\widetilde{W}=\widetilde{W}^\circ$.    
Hence, $P_{i}$ is a vertex of  $\widetilde{W}$ for any $i$ $(1\leq i \leq 2m+1)$. \\ 
\par 
\noindent
{\it Proof of  the \lq\lq if\rq\rq\; part}: \quad 
Since $P_i$ is a vertex of $\widetilde{W}$, we have the following:   
\[
\widetilde{W}=
\left\{\left.\frac{\sum_{i=1}^{k}t_{i}P_{i}}{||\sum_{i=1}^{k}t_{i}P_{i}||}\; 
\right|\; P_{i}\in X, \sum_{i=1}^{k}t_{i}=1, t_{i}\geq 0\right\}.
\]  
Thus, by Maehara's Lemma, we have the following:  
\[
\widetilde{W}^\circ=
\left\{\left.\frac{\sum_{i=1}^{k}t_{i}P_{i}}{||\sum_{i=1}^{k}t_{i}P_{i}||}\;
\right|\; P_{i}\in X, \sum_{i=1}^{k}t_{i}=1, t_{i}\geq 0\right\}^\circ
=\bigcap_{i=1}^{k}H(P_i)=\widetilde{W}.
\]
Therefore, $W$ is a self-dual Wulff shape.   
\hfill $\square$
\par
\noindent
\subsection{When is the dual Wulff shape congruent to the original Wulff shape ?}
Finally, as a generalized problem of characterization of self-dual Wulff shapes, we pose the following:   
\begin{problem}
Under what conditions is the dual Wulff shape merely congruent to the original Wulff shape  ?
\end{problem}   
We have partial results to this problem as follows:      
\begin{example}
Let $X_{2m}$ be a regular polygon with $2m$ vertices in the plane where $m\ge 2$.    
Denote the half of the length of its diagonal by $a_{2m}$.    
Suppose that the center of $X_{2m}$ is the origin and 
$a_{2m}$ satisfies the following equation:  
\[
\sin\left(\frac{\pi-\frac{2\pi}{2m}}{2}\right)=\frac{\frac{1}{a_{{}_{2m}}}}{a_{{}_{2m}}}.
\leqno{(*)}
\]
Then, $X_{2m}\ne \mathcal{D}X_{2m}$ but $\mathcal{D}X_{2m}$ is congruent to $X_{2m}$.   
\end{example}
For instance, 
consider a square $P_{1}P_{2}P_{3}P_{4}\subset \mathbb{R}^{2}$ such that  
the origin is its center and the length of its edge is $\frac{2}{a_4}$, where $a_4^{2}=\sqrt{2}$. 
Let $Q_{1}Q_{2}Q_{3}Q_{4}\subset \mathbb{R}^{2}$ be the dual Wulff shape of 
$P_{1}P_{2}P_{3}P_{4}$.   Then, $P_{1}P_{2}P_{3}P_{4}\ne Q_{1}Q_{2}Q_{3}Q_{4}$   
(see Figure \ref{figure 5}).     
And, it is easy to see that 
$Q_{1}Q_{2}Q_{3}Q_{4}$ is also a square with properties that 
the origin is its center and the length of its edge is 
$\frac{2}{a_4}$.
Thus, $Q_{1}Q_{2}Q_{3}Q_{4}$ 
is congruent to $P_{1}P_{2}P_{3}P_{4}$.    
\par 
It is not difficult to obtain the equation $(*)$ for $a_{2m}$ of general $2m$-gon $X_{2m}$.     

\begin{figure}
\begin{center}
\includegraphics[width=6cm]{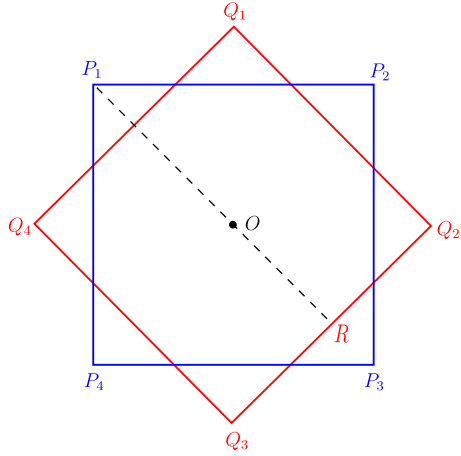}
\caption{
{Square $P_{1}P_{2}P_{3}P_{4}$ and its dual square $Q_{1}Q_{2}Q_{3}Q_{4}$. 
$P_{1}O=a_4$,} $RO=\frac{1}{a_4}$, where $a_4^{2}=\sqrt{2}$. 
}
\label{figure 5}
\end{center}
\end{figure}
}










\end{document}